\documentclass[runningheads]{llncs}
\usepackage{latexsym,amssymb,amsmath}
\renewcommand{\bbbf}{\mathbb F}

\renewcommand{\bbbz}{\mathbb Z}
\newcommand{\End}{\operatorname{End}}
\newcommand{\GL}{\operatorname{GL}}
\newcommand{\PGL}{\operatorname{PGL}}

\newtheorem{thm}{Theorem}[section]
\newtheorem{cor}[thm]{Corollary}
\newtheorem{lem}[thm]{Lemma}
\newtheorem{prop}[thm]{Proposition}

\begin{document}

\title{Some genus 3 curves with many points
\thanks{It is a pleasure to thank Hendrik Lenstra for his interest in 
this work, and for his remarks which led to Section~\ref{charLegendre} 
of this paper.}}
\titlerunning{Genus 3 curves}
\author{Roland Auer\inst{1} \and Jaap Top\inst{2}}
\institute{Department of Mathematics and Statistics, 
University of Saskatchewan,\\ 
106 Wiggins Road, Saskatoon, S7N 5E6, Canada \\ 
\email{auer@snoopy.usask.ca} \and  
IWI, Rijksuniversiteit Groningen, \\ 
Postbus 800, NL-9700 AV Groningen, The Netherlands \\ 
\email{top@math.rug.nl}}
\maketitle
\begin{abstract}
We explain a naive approach towards the problem
of finding genus 3 curves $C$ over any given finite field
${\bbbf}_q$ of odd characteristic, with a number of rational
points close to the Hasse-Weil-Serre upper bound 
$q+1+3[2\sqrt{q}]$. The method turns out to be
successful at least in characteristic $3$.
\end{abstract}

\section{Introduction}
\subsection{} The maximal number of rational points that a
(smooth, geometrically irreducible) curve of genus $g$ over
a finite field $\bbbf_q$ can have, is denoted by $N_q(g)$.
One
has the estimate (see \cite{Se})
$$N_q(g)\leq q+1+g[2\sqrt{q}]$$
in which the notation $[r]$ for $r\in\bbbr$ means the largest
integer $\leq r$. The upper bound here is called the Hasse-Weil-Serre
bound.

For $g=1$, it is a classical result of Deuring \cite{Deuring}, \cite{Wat} that
$N_q(1)=q+1+[2\sqrt{q}]$, except when $q=p^n$ with $p$ prime 
and $n\geq 3$ odd and $p$ divides $[2\sqrt{q}]$, in which case
$N_q(1)=q+[2\sqrt{q}]$. For $g=2$ an explicit formula
is due to J-P.~Serre. He stated and proved the result during a course 
\cite{SeH}
he gave at Harvard university in 1985; a nice survey including some
modifications of the original proof can be
found in Chapter~5 of the thesis \cite{Shabat}. The final result
is that if $q$ is a square and $q\neq 4,9$ then $N_q(2)= q+1+2[2\sqrt{q}]$.
Moreover $N_9(2)=20=9+1+2[2\sqrt{9}]-2$ and $N_4(2)=10=4+1+2[2\sqrt{4}]-3$. 
In case $q$ is not a square, then also $N_q(2)= q+1+2[2\sqrt{q}]$ except
when either $\gcd(q,[2\sqrt{q}])>1$ or $q$ can be written in one
of the forms $n^2+1$, $n^2+n+1$ or $n^2+n+2$. In these remaining
cases, one has that if $2\sqrt{q}-[2\sqrt{q}]\geq\frac{\sqrt{5}-1}{2}$
then $N_q(2)=q+2[2\sqrt{q}]$ and if 
$2\sqrt{q}-[2\sqrt{q}]<\frac{\sqrt{5}-1}{2}$ then $N_q(2)=q+2[2\sqrt{q}]-1$.

For $g\geq 3$ no such result is known. The best known lower bounds in
case $g\leq 50$ and $q$ a power of $2$ or $3$ which is $\leq 128$ can
be found in \cite{G-V}. In \cite[\S~4]{SerreBordeaux} J-P.~Serre gives values 
of
$N_q(3)$ for $q\leq 19$ and for $q=25$. Moreover he shows in
\cite[p.~64-69]{SeH}  that $N_{23}= 48$. Hence we have the
following table.
$$\begin{array}{|c|c|c|c|c|c|c|c|c|c|c|c|c|c|c|c|c|c|c|c|c|}\hline
q&2&3&4&5&7&8&9&11&13&16&17&19&23&25&27&29&31&32&37&41\\\hline
N_q(3)&7&10&14&16&20&24&28&28&32&38&40&44&48&56&56&60&\geq 56&64&\geq 68&\geq 
72\\\hline
\end{array}
$$
The entries for $q=29,31,37$ are obtained using the technique from 
the current paper;
its main goal is to
give lower bounds for $N_q(3)$ by restricting ourselves
to one specific family of curves of genus $3$.

\subsection{} Let $k$ be a field of characteristic different from $2$.
The plane quartic $C_{\lambda}$ given by
$x^4+y^4+z^4=(\lambda+1)(x^2y^2+y^2z^2+z^2x^2)$ is for $\lambda
\in k$ with $\lambda\neq -3, 1, 0$ a geometrically irreducible,
smooth curve of genus $3$. The degree $4$ polynomials given here
are fixed by the subgroup $G<\PGL(3,k)$ generated by
$\sigma: (x,y,z)\mapsto (y,z,x)$ and $\tau: (x,y,z)\mapsto(y,-x,z)$.
The group $G$ is isomorphic to $S_4$, the symmetric group
on $4$ elements. Hence $G$ is contained in the group
of automorphisms of $C_\lambda$. For general $\lambda$ the
automorphism group of $C_\lambda$ in fact equals $G$.

These curves occur in the classification of non-hyperelliptic
genus $3$ curves with nontrivial automorphism group, as given
in \cite[p.~2.88]{Henn} and in \cite[Table~5.6, pp.~63-64]{Vermeulen}.

Suppose $\lambda\neq 0,1$. By $E_\lambda$ we denote the
elliptic curve given by the equation $y^2=x(x-1)(x-\lambda)$.
If moreover $\lambda\neq -3$ then we write $E^{(\lambda+3)}_\lambda$
for the elliptic curve with equation $(\lambda+3)y^2=x(x-1)(x-\lambda)$.
The relation with the curves $C_\lambda$ is as follows.

\begin{lem}\label{jacobian} 
Suppose $k$ is a field of characteristic 
different from $2$ and $\lambda\in k\setminus\{0,1,-3\}$. Then
the jacobian of the curve $ C_{\lambda}$ given by
$x^4+y^4+z^4=(\lambda+1)(x^2y^2+y^2z^2+z^2x^2)$ is over $k$
isogenous to the product
$E^{(\lambda+3)}_\lambda\times E^{(\lambda+3)}_\lambda
\times E^{(\lambda+3)}_\lambda$, where $E^{(\lambda+3)}_\lambda$
denotes the elliptic curve with equation $(\lambda+3)y^2=x(x-1)(x-\lambda)$.
\end{lem}

\noindent
{\sl Proof.} Most of this is shown in \cite[pp.~40-41]{To}; one takes
the quotient of $C_{\lambda}$ by the involution $(x,y,z)\mapsto(-x,y,z)$.
The resulting curve has genus $1$ and it admits an involution
without any fixed points. Taking the quotient again results in
an elliptic curve, given by
$y^2=x^3+2(\lambda+1)(\lambda+3)x^2+(\lambda-1)(\lambda+3)x$.
The $2$-isogeny with kernel generated by $(0,0)$ maps this
curve onto $E_\lambda^{(\lambda+3)}$ (compare the
formulas for $2$-isogenies as given in \cite[III \S~4]{Tate-Sil}). Write
$\pi:C_{\lambda}\rightarrow E_\lambda^{(\lambda+3)}$ for
the composition of all maps described here.
Then
$$\rho=(\pi,\pi\sigma,\pi\sigma^2):\; C_{\lambda}\rightarrow
E^{(\lambda+3)}_\lambda\times E^{(\lambda+3)}_\lambda
\times E^{(\lambda+3)}_\lambda$$
where $\sigma: (x,y,z)\mapsto (y,z,x)$ is one of the
automorphisms of $C_\lambda$. The fact that $\rho$ induces
an isomorphism between the spaces of regular $1$-forms implies
that $\rho$ induces an isogeny between $\mbox{Jac}(C_{\lambda})$
and the triple product of $E^{(\lambda+3)}_\lambda$. {\hfill{$\Box$}}

\begin{cor}\label{comparetraces}
With notations as above, one finds for $\lambda\in\bbbf_q$ 
with $q$ odd and $\lambda\neq 0,1,-3$ that
$$\#C_{\lambda}(\bbbf_q)=3\#E_\lambda^{(\lambda+3)}(\bbbf_q)-2q-2.$$
\end{cor}

\noindent
{\sl Proof.} It is a well known fact that
$\# C_{\lambda}(\bbbf_q)$ equals $q+1-t$, where $t$ is the
trace of Frobenius acting on a Tate module of $\mbox{Jac}(C_{\lambda})$.
Lemma~\ref{jacobian} implies that this Tate module is isomorphic
to a direct sum of three copies of the Tate module of 
$E_\lambda^{(\lambda+3)}$.
Hence $t=3t'$ where $t'$ is the trace of Frobenius on the Tate
module of $E_\lambda^{(\lambda+3)}$. Since this trace equals
$q+1-\#E_\lambda^{(\lambda+3)}(\bbbf_q)$, the result follows.
{\hfill{$\Box$}}

\subsection{}\label{strategy}
Our strategy for finding a curve of genus $3$ over a finite field
$\bbbf_q$ with odd characteristic should now be clear: find 
$\lambda$ such that $\#E_\lambda^{(\lambda+3)}(\bbbf_q)$ is as
large as possible and use Corollary~\ref{comparetraces}.
This works quite well for small $q$, using a direct search.
In fact, as will be explained in Section~\ref{charp} below,
it is not even necessary here to calculate 
$\#E_\lambda^{(\lambda+3)}(\bbbf_q)$ for many values 
$\lambda\in\bbbf_q$.

We obtain a general result when the characteristic of
$\bbbf_q$ equals $3$, because in that case we deal with a
curve $E_\lambda^{(\lambda)}$ which is isomorphic to the
curve $E_\mu$ with $\mu=1/\lambda$. Since it is precisely
known which values $\#E_\mu(\bbbf_q)$ attains (see \cite{A-T}
and also Section~\ref{charLegendre} below), one obtains
a nice explicit lower bound for $N_{3^n}(3)$. In fact, the
result implies that the difference between $N_{3^n}(3)$ and
the Hasse-Weil-Serre bound is bounded independently of $n$:

\begin{prop}\label{resultchar3}
For every $n\geq 1$ the inequality
$$3^n+1+3[2\sqrt{3^n}]-N_{3^n}(3)\leq\left\{
\begin{array}{cl}
0 &\mbox{ if }\; n\equiv 2\bmod 4;\\
12 &\mbox{ if }\; n\equiv 0\bmod 4;\\
21 &\mbox{ if }\; n\equiv 1\bmod 2\end{array}\right.
$$ holds.
\end{prop}

\noindent
For the proof we refer to Section~\ref{char3}. Note that
this proves a special case of a conjecture of 
J-P.~Serre \cite[p.~71]{SeH},
which says that for {\em all} $q$ the difference
$q+1+3[2\sqrt{q}]-N_q(3)$ should be bounded independently of $q$.
 
\subsection{}
In characteristic at least $5$ we have not been able to obtain
a general result such as given in Proposition~\ref{resultchar3}.
However, the fact that a curve $E_\lambda^{(\lambda+3)}$ is
either isomorphic to $E_{\lambda}$ or it is a quadratic twist
of $E_{\lambda}$, implies (again using \cite{A-T}) that for
every finite field $\bbbf_q$ of odd characteristic, a curve
$C_\lambda$ as above exists for which $\#C(\bbbf_q)$ is
at most $21$ off from either the Hasse-Weil-Serre upper bound
$q+1+3[2\sqrt{q}]$, or from the analogous lower bound
$q+1-3[2\sqrt{q}]$. This is proven in Section~\ref{charp}.
We note that a sharper result of the same kind (with $21$ replaced
by $3$) was obtained by Kristin Lauter \cite{kristin}, \cite{LauS} using an 
entirely different method. 

As Everett Howe pointed out to us, it is in fact possible to improve
our result slightly by replacing the product
$ E_\lambda^{(\lambda+3)}\times E_\lambda^{(\lambda+3)}\times 
E_\lambda^{(\lambda+3)}$
we use, by a product $E\times E\times E_{\lambda}$ in which
$E/\bbbf_q$ is an elliptic curve with a rational point of order $2$ and
$\#E(\bbbf_q)$ maximal under that condition, and $E_{\lambda}/\bbbf_q$ is
a Legendre elliptic curve over $\bbbf_q$ with as many rational points
as possible. The result of Everett Howe, Franck Lepr\'{e}vost and
Bjorn Poonen \cite[Prop.~15]{HLP} in this case implies that
either this product or its standard quadratic twist is isogenous over
$\bbbf_q$ to the jacobian of a smooth genus $3$ curve over $\bbbf_q$.
It may be noted that the estimate obtained in this way is in general still 
weaker
than Lauter's result (it replaces our $21$ by $9$ instead of by $3$).

\section{A characterization of Legendre elliptic curves}\label{charLegendre}
Suppose $K$ is a field of characteristic $\neq 2$, and $E/K$ is
an elliptic curve. We will say that $E/K$ is a Legendre elliptic curve
over $K$ if there is a $\lambda\neq 0,1$ in $K$ such that $E$ is over $K$
isomorphic to $E_{\lambda}$ given by $y^2=x(x-1)(x-\lambda)$.
A necessary but in general not sufficient condition for an
elliptic curve $E/K$ to be a Legendre elliptic curve over $K$ is that
all points of order $2$ on $E$ are $K$-rational.
An  intrinsic description of Legendre elliptic curves is given
as follows. Take a separable closure $K^{\mbox{\scriptsize sep}}$ of $K$
and write $G_K=\mbox{Gal}(K^{\mbox{\scriptsize sep}}/K)$ for its
Galois group.

\begin{lem}\label{Legendrecriterion}
The statements
\begin{enumerate}
\item $E$ is a Legendre elliptic curve over $K$;
\item $E$ can be given by an equation $y^2=(x-a)(x-b)(x-c)$ in which
at least one of $\pm(a-b),\pm(b-c),\pm(c-a)$ is a square in $K^*$;
\item $E$ has all its points of order $2$ rational over $K$, and there
exists a point $P\in E(K^{\mbox{\scriptsize sep}})[4]$ such that
$-P$ is not in the $G_K$-orbit of $P$.
\end{enumerate}
are equivalent.
\end{lem}

\noindent 
{\sl Proof}. The equivalence of (1) and (2) is easy. To verify that (2)
and (3) are equivalent, suppose (after possibly permuting $a,b,c$) that
$a-b$ is a square and that $E$ is given by $y^2=(x-a)(x-b)(x-c)$.
The point $T_b=(b,0)$ in $E(K)$ has order $2$, and the quotient 
$E':=E/\langle T_b\rangle$ admits an isogeny of degree $2$:
$\varphi:E'\rightarrow E$ defined over $K$ (the dual isogeny of the
quotient map). A very well known property 
(compare \cite[III \S~5]{Tate-Sil})
of $\varphi$ is that
the image $\varphi(E'(K))\subset E(K)$ equals the kernel of the
homomorphism $E(K)\rightarrow K^*/{K^*}^2$ defined by $T_b\mapsto
(b-a)(b-c)$ and $(x,y)\mapsto x-b$ for all $(x,y)\in E(K)$ with
$(x,y)\neq T_b$. Hence the condition that $a-b$ be a square is
equivalent with the property that the point $T_a:=(a,0)\in E(K)$
is in the image of $E'(K)$. This means precisely that a pair
of points $\{P,P+T_b\}\subset E$ exists which is $G_K$-stable, and
$2P=T_a$.
Hence $P$ is a point of order $4$ on $E$, and for all $\sigma\in G_K$ 
we have $\sigma(P)-P\in \{O,T_b\}$. In particular $\sigma(P)-P\neq 2P$,
which means $\sigma(P)\neq -P$ for all $\sigma\in G_K$.

Vice versa, suppose given a point $P$ of order $4$ with the property
$\sigma(P)\neq -P$ for all $\sigma\in G_K$. 
Since all $2$-torsion of $E$ is $K$-rational, we have that
$\sigma(P)-P\in E(K)[2]$ and moreover the condition $\sigma(P)\neq -P$
implies that $\sigma(P)-P$ is in a cyclic subgroup of $E(K)[2]$ which
is independent of $\sigma$. Hence we have points $T$ and $2P$ of order $2$,
where $\{P,P+T\}$ is $G_K$-stable. As we have seen, this implies the
statements (1) and (2).
\hfill{$\Box$}

\begin{cor}\label{LegFrob}
Suppose $q$ is a power of an odd prime and $E/{\bbbf}_q$ is an
elliptic curve. Let $\pi\in \End(E)$ be the Frobenius endomorphism
(raising coordinates to the power $q$).
Then $E$ is a Legendre elliptic curve over ${\bbbf}_q$ if and only if
$\pi+1\in 2\End(E)$ but $\pi+1\not\in 4\End(E)$.
\end{cor}

\noindent 
{\sl Proof}. 
The Galois group $G_{\bbbf_q}$ is topologically generated by
the $q$-th power map, and this generator acts on $E$ via the
endomorphism $\pi$.
The condition $\pi+1\in 2\End(E)$ is equivalent with the
statement that $E$ has all its points of order $2$ rational over $\bbbf_q$.
In the same manner, the condition $\pi+1\not\in 4\End(E)$
precisely means that a point $P$ of order $4$ exists, with the
property $\pi(P)\neq -P$. 
Since the Galois group $G_{\bbbf_q}$ acts on $E(\overline{\bbbf}_q)[4]$ 
via a (cyclic) subgroup of the kernel of
$\GL_2(\bbbz/4\bbbz){\stackrel{\mbox{\scriptsize mod 2}}
{\longrightarrow}}\GL_2(\bbbz/2\bbbz)$, it follows that
$\sigma(P)\neq -P$ for all $\sigma\in G_{\bbbf_q}$.
Hence Lemma~\ref{Legendrecriterion} implies that $E$ is a Legendre
elliptic curve over $\bbbf_q$.

Vice versa, if $E$ is a Legendre
elliptic curve over $\bbbf_q$, then by Lemma~\ref{Legendrecriterion}
we know that $P\in E(\overline{\bbbf_q})[4]$ exists with $\pi(P)\neq -P$,
which implies that $\pi+1$ is not divisible by $4$ in $\End(E)$.
We have that $\pi+1\in 2\End(E)$ since $\pi$ acts trivially on
all points of order $2$.

This proves the corollary.\hfill{$\Box$}

\begin{prop}\label{ATresult}
An elliptic curve $E/\bbbf_q$ (with $q$ odd) for which $\#E(\bbbf_q)\in
4\bbbz$ is isogenous to a Legendre elliptic curve over $\bbbf_q$, except
in the following case: $q=r^2$ with $r\in 1+4\bbbz$, and 
$\#E(\bbbf_q)=q+1+2r$.
\end{prop}

\noindent 
{\sl Proof}. 
(This result was first presented in \cite{A-T}, however, with a somewhat 
different 
proof. The present proof is more conceptual, but it gives less information
concerning the possible values of Legendre parameters $\lambda$ in the 
supersingular case.)

Let $\pi\in \End(E)$ be the ($q$-th power) Frobenius. The proof
considers two cases.

First, suppose $\pi=r\in \bbbz$. Then $q=\deg(\pi)=r^2$ and
$\#E(\bbbf_q)=(r-1)^2$. Any curve $E'$ isogenous to $E$ then also
satisfies $\#E'(\bbbf_q)=(r-1)^2$ and Frobenius in $\End(E')$ is
equal to $r$. By Corollary~\ref{LegFrob}, one (and equivalently, all of them)
such curve $E'$ is a
Legendre elliptic curve over $\bbbf_q$ precisely when $r+1$ is even, but not
divisible by $4$. The latter condition is equivalent with $r\equiv 1\bmod 4$.
This proves the statement in the case $\pi\in\bbbz$.

If $\pi\not\in\bbbz$ then $\bbbz[\pi]\subset\End(E)$ is an order
in the ring of integers of an imaginary quadratic field $K$. We have that
$\#E(\bbbf_q)=(1-\pi)(\overline{1-\pi})$ where the bar denotes complex
conjugation in $K$. The condition $\#E(\bbbf_q)\equiv 0\bmod 4$ implies
that $(1-\pi)/2$ is integral. Now consider the order $A:=\bbbz[(1+\pi)/2]$.
By construction, $\pi\in A$ satisfies $\pi+1\in 2A$ and $\pi+1\not\in 4A$.
It is a result of Waterhouse \cite[Thm.~4.5]{Wat} (compare \cite[p.~194]{Sch}
where a mistake in the original result is corrected), that a curve 
$E'/\bbbf_q$
exists with an isomorphism $\End(E')\cong A$ such that under this 
isomorphism
Frobenius on $E'$ corresponds to $\pi\in A$. This implies in
particular that $\#E'(\bbbf_q)=\#E(\bbbf_q)$ and hence $E'$ and $E$
are isogenous. Moreover, using Corollary~\ref{LegFrob} we know
that $E'$ is a Legendre elliptic curve over $\bbbf_q$.
This proves the proposition.\hfill{$\Box$}

\section{Characteristic $3$}\label{char3}
We will now prove Proposition~\ref{resultchar3}. Take $n\geq 1$ and write
$q:=3^n$, $m:=[2\sqrt{q}]$ and $q+1+m=N+r$ with $N\in 4\bbbz$ and $0\leq 
r\leq 3$.
As explained in (\ref{strategy}), we will examine how close to the
upper bound $q+1+m$ the number of $\bbbf_q$-points on a Legendre elliptic
curve $E_{1/\lambda}\cong E^{(\lambda)}_{\lambda}$ can be, for 
$\lambda\in\bbbf_q$.

If $n$ is odd and moreover $N\equiv 1\bmod 3$ (the smallest $n$
where this is the case, is $n=11$ which gives $m=841$ and
$N=3^{11}+1+840$), then we replace $N$ by $N-4$.
The resulting number $N$ satisfies $q+1-m\leq N\leq q+1+m$, 
and moreover we know from \cite{Deuring}
that $E/\bbbf_q$ exists with $\#E(\bbbf_q)=N$. If $n$ is odd, then
Proposition~\ref{ATresult} implies the existence of $\lambda\in\bbbf_q$
with $\#E^{(\lambda)}_{\lambda}(\bbbf_q)=N$. Hence 
Corollary~\ref{comparetraces}
yields a genus $3$ curve $C_{\lambda}$ with $\#C_{\lambda}(\bbbf_q)=3N-2q-2$.
In particular, this shows that 
$$q+1+3m-N_q(3)\leq q+1+3m-3N+2q+2=3r+12\leq 21$$
for odd $n$ (in fact, even $\leq 3r\leq 9$ unless $m$ is divisible by $3$).

If $n$ is even, then $m=2\cdot 3^{n/2}$ and (again using Deuring's results
\cite{Deuring}) an elliptic curve $E/\bbbf_q$ exists with 
$\#E(\bbbf_q)=q+1+m$.
By Proposition~\ref{ATresult}, this number of points occurs for a
Legendre elliptic curve only in case $m/2\equiv 3\bmod 4$, i.e., when
$n\equiv 2\bmod 4$. Hence under this condition we obtain a curve
$C_{\lambda}$ whose number of points attains the Hasse-Weil-Serre bound.

In the remaining case we have $n\equiv 0\bmod 4$. Here the number $q+1+m$
does not occur as $\#E_{\lambda}^(\lambda)(\bbbf_q)$, for any 
$\lambda\in\bbbf_q$.
Hence we take the largest smaller possibility, which is
$q+1+m-4$. Proposition~\ref{ATresult} implies that a Legendre elliptic
curve with this number of points over $\bbbf_q$ indeed occurs. It
follows that a genus $3$ curve $C_{\lambda}/\bbbf_q$ exists
with $\#C_{\lambda}(\bbbf_q)=3(q+1+m-4)-2q-2=q+1+3m-12$.
This implies the inequality given in Proposition~\ref{resultchar3}.
\hfill{$\Box$}

\section{Examples in characteristic $>3$}\label{charp}
The problem which arises when one attempts to adapt the argument presented
in Section~\ref{char3} to finite fields of characteristic $>3$, can
already be seen in the following result.

\begin{prop}
Suppose $q$ is a power of a prime $p>3$, and $m:=[2\sqrt{q}]$. 
Over $\bbbf_q$, a curve
$C_{\lambda}$ of genus $3$ exists such that either
$\#C_{\lambda}(\bbbf_q)\geq q+1+3m-21$ or
$\#C_{\lambda}(\bbbf_q)\leq q+1-3m+21$.
\end{prop}

\noindent
As we mentioned in the introduction, a somewhat stronger result has been
obtained by Kristin Lauter \cite{kristin}, \cite{LauS} using quite different 
techniques.
Moreover a variant of our proof may be obtained by using a result of
Everett Howe, Franck Lepr\'{e}vost and Bjorn Poonen
\cite[Prop.~15]{HLP}.

\noindent
{\sl Proof.} Write $q+1+m=N+r$ with $N\in 4\bbbz$ and $0\leq r\leq 3$.
Then one of $N,2q+2-N$ occurs as the number of points on a Legendre
elliptic curve $E_{\lambda}/\bbbf_q$, except possibly when $r>0$ and
$p$ divides $m-r$. In that case, we replace $N$ by $N':=N-4$ and we
obtain a number of points which does occur. 

This gives us an elliptic curve $E_{\lambda}$. The corresponding curve
$E^{(\lambda+3)}_{\lambda}$ has either $N$ or $N'$ points, or in case
$\lambda+3$ is not a square in $\bbbf_q$ this number is $2q+2-N$ or
$2q+2-N'$. Since this number is at distance at most $7$ from one
of $q+1\pm m$, Corollary~\ref{comparetraces} implies the result.
\hfill{$\Box$}

\begin{prop}
Suppose $p\equiv 3\bmod 4$ is a prime number, $n\geq 1$ is an odd
integer and $q=p^{2n}$.
Then $N_q(3)=q+1+6p^n$ equals the Hasse-Weil-Serre bound.
\end{prop}

\noindent
{\sl Proof.} Take $\lambda=-1\in\bbbf_p$. Since $p\equiv 3\bmod 4$,
 the elliptic curve $E_{\lambda}/\bbbf_p$ is supersingular.
This implies $\#E_{\lambda}(\bbbf_p)=p+1$ (in case $p=3$, this follows from
the fact that the number of points is a multiple of $4$, and also
of course from a direct calculation). One concludes 
that $\#E_{\lambda}(\bbbf_q)=q+1+2p^n$. Since $\lambda+3\neq 0$ as an
element of $\bbbf_p$ is a square in $\bbbf_q$, the two curves
$E_{\lambda}$ and $E_{\lambda}^{(\lambda)}$ are isomorphic over $\bbbf_q$.
Corollary~\ref{comparetraces} therefore yields that the genus $3$
curve $C_{\lambda}$ attains the Hasse-Weil-Serre bound over $\bbbf_q$.
\hfill{$\Box$}

\vspace{\baselineskip}
Note that the genus $3$ curve used in the above proposition is in fact
the famous Fermat quartic. Hence the result is probably well known.

\subsection{}
In practice, a fairly efficient method to find $\lambda\in\bbbf_q$ for which
$\# E^{(\lambda+3)}_{\lambda}(\bbbf_q)$ equals a given number
$N\equiv 0\bmod 4$ can be given in case $q=p$ a prime or $q=p^2$ the square
of a prime. This works as follows. Write $N=q+1-t$. 

We first treat the case $q=p^2$ and  $t=\pm 2p$. Exactly
one of the two numbers $p^2+1\pm 2p$ occurs as a number of points
of a Legendre elliptic curve over $\bbbf_{p^2}$, and this number
is attained in our family precisely for the supersingular $\lambda\neq -3$
such that $\lambda+3$ is a square in $\bbbf_{p^2}$; the number with the 
opposite
choice of sign occurs for the ones such that $\lambda+3$ is a nonsquare.

In the remaining cases, the Hasse inequality tells us $|t|< 2\sqrt{q}$. 
Hence we find
$t$ exactly if we know $t\bmod 4p$. Now $t\bmod 4$ is already known, hence
it suffices to find a $\lambda$ such that 
$\# E^{(\lambda+3)}_{\lambda}(\bbbf_q)\equiv q+1-t\bmod p$.
If we write $\chi:\bbbf_q^*\rightarrow \pm 1$ for the nontrivial
character with kernel ${\bbbf_q^*}^2$, this means we look for $\lambda\neq -3$
such that $\#E_{\lambda}(\bbbf_q)\equiv 1-\chi(\lambda+3)t$.
It is well known \cite[V \S~4]{Sil} that 
$\#E_{\lambda}(\bbbf_q)\equiv 1-\left((-1)^{(p-1)/2}H_p(\lambda)\right)^e$,
with $e=1$ if $q=p$ and $e=1+p$ if $q=p^2$. Here
$H_p(\lambda)=\sum_{i=0}^{(p-1)/2}  \binom{(p-1)/2)}{i}^2\lambda^i$
is the so-called Hasse polynomial, whose coefficients can be
computed using an easy recursion.
Hence in case $q=p$ we have to solve $H_p(\lambda)=\pm t$ for 
$\lambda\in\bbbf_p$, and
then check whether $\chi(\lambda+3)$ has the correct value.
Similarly, when $q=p^2$ we look for solutions in $\bbbf_q$ of
$H_p(\lambda)H_p(\lambda^p)=\pm t$. This works reasonably efficient
for -say- $q< 10^7$.

\subsection{}
Using the package KANT, we tested which values $\#E(\bbbf_q)\equiv 0\bmod 4$
occur for the curves $ E^{(\lambda+3)}_{\lambda}/\bbbf_q$, for all odd
$q<100000$. It turns out that for most $q$, all values are attained. In the 
table
below, we list all $q<100000$ where this is {\em not} the case, and for
each of them the missing value(s) $\#E(\bbbf_q)$. As can be seen from the
data, usually there is only one such missing value, which moreover is
always one of the minimal or the maximal possible number $\equiv 0\bmod 4$
in the interval $\left[q+1-[2\sqrt{q}],q+1+[2\sqrt{q}]\right]$. We list a 
sign $\pm$
indicating which of these possibilities occurs for a missing value.
There are exactly two exceptional cases for $q<100000$.
The first one is $q=7^4$: here two values don't occur,
namely the maximum $2500=7^4+1+2\cdot 7^2$ and also $2396=7^4+1-6$.
The other one is $q=5^6$. The two values missing here are
the minimal  one $5^6+1- 2\cdot 5^3$ and also $15380=q+1-2\sqrt{q}+4$.
The following table gives all other $q<100000$ with a missing value.
$$\begin{array}{|c|c||c|c||c|c||c|c||c|c||c|c||c|c||c|c|}
q&\pm&q&\pm&q&\pm&q&\pm&q&\pm&q&\pm&q&\pm&q&\pm\\\hline
5&+&7&-&3^2&-&13&-&19&-&5^2&-&7^2&-&67&1\\
3^4&+&5^3&-&13^2&-&173&+&293&+&7^3&-&487&-&23^2&-\\
5^4&+&3^6&-&733&-&787&+&907&+&2503&+&3253&+&4493&-\\
4903&-&5333&+&5479&-&5779&-&3^8&+&7573&-&9413&+&10639&-\\
11239&-&11243&+&12547&-&11^4&+&14887&-&17959&+&18773&+&23719&+\\
24967&-&25603&-&27893&-&13^4&+&31687&-&33287&+&33493&-&37253&+\\
42853&-&46663&-&51991&+&52903&-&58567&+&3^{10}&-&64013&+&65539&+\\
67607&+&71293&-&76733&+&17^4&+&85853&+&92419&+&94253&-&99859&-\\\hline\end{array}
$$
The table shows that for all but $30$ values $q<100000$, the maximal
value of $\#C_{\lambda}(\bbbf_q)$ equals $q+1+3t$, where $q+1+t\equiv 0\bmod 
4$
is the maximal number of points of an elliptic curve over $\bbbf_q$ with all 
its
points of order $2$ rational.

Whenever the Hasse-Weil-Serre bound is divisible by $4$, we may be in
the lucky circumstance that it is reached using our family of
curves. This happens quite frequently, for instance when $q$ equals any 
of the primes
$19,\;29,\;53$, $67,\;71,\;89,\;103$,
$107,\;151,\ldots$. In the case $q=173$ the bound $q+1+3m$
is a multiple of $4$, but as the table above shows, our curves do not
attain it.

The data seems to indicate that for much more than $50\%$ of the prime powers 
$q$,
all possible values $N_q\equiv 0\bmod 4$ are attained by the family
$E_{\lambda}^{(\lambda+3)}$. Moreover, the only occurrences of a $q$ for 
which more than 
one value is missing, happened at `high' even powers of a prime number.
We have no theoretical explanation for this. A numerical test over all 
$q=p^{2n}<10^7$
revealed exactly one more case where two values are missing, namely
at $q=7^6$. We have also not been able to explain why in all
cases where we found that exactly one value is missing, this missing value is 
one
of the maximal or minimal number of points.

\end{document}